\documentclass[a4paper,12pt,leqno]{article}
\usepackage{latexsym}
\usepackage[all]{xy}

\usepackage{amssymb} 
\usepackage{amsmath} 

\usepackage{tikz}
\usetikzlibrary{arrows,cd}  
\usepackage{amsthm}

\def\K{{\mathbb{K}}}
\def\R{{\mathbb{R}}}
\def\C{{\mathbb{C}}}
\def\A{{\mathcal{A}}}
\def\B{{\mathcal{B}}}
\def\C{{\mathcal{C}}}

\DeclareMathOperator{\Der}{Der}

\DeclareMathOperator{\Shi}{Shi}
\DeclareMathOperator{\Cat}{Cat}

\DeclareMathOperator{\POexp}{POexp}

\numberwithin{equation}{section}

\newcommand{\owari}{\hfill$\square$}

\newtheorem{theorem}{Theorem}[section]
\newtheorem{prop}[theorem]{Proposition}
\newtheorem{cor}[theorem]{Corollary}

\newtheorem{define}[theorem]{Definition}

\theoremstyle{remark}

\newtheorem{example}[theorem]{Example}

\title{Free paths of arrangements of hyperplanes}
\author{Takuro Abe\footnote{
Department of Mathematics,
Rikkyo University,
3-34-1, Nishi-Ikebukuro,
Toshima-ku,
1718501, Tokyo.
Email:abetaku@rikkyo.ac.jp\ \ 
ORCID:0000-0002-4477-8450} and Toru Yamaguchi\footnote{
RICOH COMPANY, LTD., 
Nakamagome 1-3-6 Ota-ku, Tokyo
Email:Toru.1.Yamaguchi@jp.ricoh.com
}}
\date{\today}

\begin{document}

\maketitle

\begin{abstract}
We study the free path problem, i.e., if we are given two free arrangements of hyperplanes, then we can connect them by free arrangements or not. We prove that if an arrangement $\A$ and 
$\A \setminus \{H,L\}$ are free, then at least one of two among them is free. When $\A$ is in the three dimensional arrangement, we show a stronger statement.
\end{abstract}

\section{Introduction}
Let $\K$ be a field, $V=\K^\ell$, $S:=\mbox{Sym}^*{V^*} =\K[x_1,\ldots,x_\ell]$ for a basis $x_1,\ldots,x_\ell$ for $V^*$ and let 
$\Der S:=\oplus_{i=1}^\ell S \partial_{x_i}$. Let $\A$ be an \textbf{arrangements of hyperplanes}, i.e., a finite set of linear hyperplanes in $V$. For each $H \in \A$ let us fix a defining linear form $\alpha_H \in V^*$ and let $Q(\A):=\prod_{H \in \A} \alpha_H$. Then the 
\textbf{logarithmic derivation module} $D(\A)$ of $\A$ is defined by 
$$
D(\A):=\{\theta \in \Der S \mid \theta(\alpha_H) \in S \alpha_H\ (\forall H \in \A)\}. 
$$
$D(\A)$ is an $S$-graded reflexive module of rank $\ell$, but not free in general. So we say that $\A$ is 
\textbf{free} with \textbf{exponents} $\exp(\A)=(d_1,d_2,\ldots,d_\ell)$ if 
$$
D(\A)\simeq \oplus_{i=1}^\ell S[-d_i].
$$
In this article we assume that $\A$ is essential, i.e., $\cap_{H \in \A} H =\{0\}$.  
In this case it is known that $$
D(\A) =S \theta_E \oplus D_H(\A)
$$
for any $H \in \A$, where 
$$
D_H(\A):=\{\theta \in D(\A) \mid \theta(\alpha_H)=0\}.
$$
So $1$ is always contained in $\exp(\A)$ that is the degree of the Euler derivation 
$\theta_E$. 

In the study of hyperplane arrangements, free arrangements 
have played the key role which connects algebra, topology, combinatorics and geometry of hyperplane arrangements. In particular, how to find and construct a free arrangement is very fundamental and important. Among them, the most 
famous result is Terao's addition-deletion theorem (Theorem \ref{adddel}). This describes the freeness of $\A$ and $\A':=\A \setminus \{H\}$ combined with the data of $\A^H:=\{L\cap H \mid L \in \A'\}$ that is an arrangement in $H \simeq \K^{\ell-1}$. On the other hand, in \cite{AT2}, the concept of the \textbf{free filtration} was introduced, i.e., for a free arrangement $\A$, we say that 
the filtration 
$$
\emptyset=\A_0 \subset \A_1 \subset \cdots \subset \A_n=\A$$
is the free filtration if $|\A_i \setminus \A_{i-1}|=1$ for all $i$, and 
$\A_i$ is free for all $i$. 
So a free arrangement with a free filtration is the arrangement which can be constructed by using the addition theorem (Theorem \ref{adddel}) from empty arrangements. To such an arrangement, it was proved in \cite{A6} that Terao's conjecture asserting that the freeness is combinatorially determined is true (see Theorem 
\ref{addcombin}). Here the combinatorics 
of arrangements is the poset structure of the \textbf{intersection lattice} defined by 
$$
L(\A):=
\{\cap_{H \in \B} H \mid 
\B \subset \A\}.
$$
Then more explicitly Theorem \ref{addcombin} asserts that if the freeness of 
$\B$ is free and it is determined by $L(\B)$, and there is a ``free filtration'' from $\B$ to $\A$, then the freeness of $\A$ is also determined by $L(\A)$. 
In these days, combinatorial dependency of many class of free hyperplane arrangements have been proved 
like divisionally, additively, and stair-free arrangement, see 
\cite{A2} and \cite{A6} for details. 
So now it is important to study whether we can connect 
two free arrangements by free arrangements among them. 
So let us introduce a generalized concept of this free filtration as follows:

\begin{define}
Let $\B \subset \A$ be both free. Then we say that there is a \textbf{free path} from $\B$ to $\A$ if there is a filtration 
$$
\B=
\A_0 \subset \A_1 \subset \cdots \subset \A_n=\A$$
such that 
$|\A_i \setminus \A_{i-1}|=1$ for all $i$, and 
$\A_i$ is free for all $i$.
\label{freepathdef}
\end{define}

So we want to ask for given two free arrangements, we can connect them by free arrangements or not. 
Since it is not easy to investigate the structure of logarithmic derivation modules of two arrangements 
$\A \supset \B$ such that $|\A \setminus \B|>1$, there have been no such research. In this article, first we prove the following fundamental result on the free path:

\begin{theorem}
Let $\A \supset \{H_1,H_2\},\ 
\A_i:=\A \setminus \{H_i\}\ (i=1,2)$ and let 
$\B:=\A \setminus \{H_1,H_2\}$. If $\A$ and $\B$ are both free, 
then at least one of $\A_1$ and $\A_2$ is free.
\label{freepath}
\end{theorem}

So if $|\A \setminus \B|=2$, then there is always a free path. Also if $\ell=3$, then we can extend this result as follows:

\begin{theorem}
Let $\A$ be free in $\K^3$, $H_1,H_2,H_3 \in \A$ be 
distinct hyperplanes. Let $\A_i:=\A \setminus \{H_i\}$,\ 
$\B:=\A \setminus \{H_1,H_2,H_3\}$ and $\B_i:=\B \cup 
\{H_i\}$. Assume that $\B$ is free too.
Then there is a 
free path between $\A$ and $\B. $
\label{main3}
\end{theorem}

Theorem \ref{main3} is sharp, i.e., if $|\A \setminus \B| \ge 4$, then we can find 
a counterexample to Theorem \ref{main3}, see Example \ref{pentagon}. 

The proof of Theorems \ref{freepath} and \ref{main3} is based on the SPOG theory, which is a good structure 
of the logarithmic derivation modules close to the freeness. See Definition \ref{SPOGdefine}
and Theorems \ref{SPOG}, \ref{SPOG2} for the SPOGness. By using them, we can 
connect $\A \supset \B$ such that $2 \le |\A \setminus \B | \le 3$. 

The organization of this article is 
as follows. In \S2 we introduce several definitions and results to prove our main theorem. In \S3 we prove the main results, give some applications and examples.
\medskip

\noindent
\textbf{Acknowledgements}. The authors are grateful to Lukas 
K\"{u}hne for suggesting the ``free path'' problem, and 
several fruitful discussions. 
 The first author 
is partially supported by JSPS KAKENHI Grant Number JP21H00975.

\section{Preliminaries}

Now let us gather several definitions and results for the proof of main results. A general reference in this section is \cite{OT}. 
Before results let us define a notation. For the multiset $(d_1,\ldots,d_s)$ of integers, $(d_1,\ldots,d_s)_\le$ denotes the multiset 
satisfying $d_1 \le \cdots \le d_s$. In the same way, 
$(d_1,\ldots,d_s)_<$ denotes the multiset 
satisfying $d_1 < \cdots < d_s$.

First let us recall the basic result connecting the algebraic structure of $\A$ and $\A':=\A \setminus\{H\}$. 

\begin{theorem}[\cite{T1}]
Let $H \in \A,\ \A':=\A \setminus \{H\}$. Then 
there is a polynomial $B$ of 
degree $|\A'|-
|\A^H|$ such that 
$$
\theta(\alpha_H) \in (\alpha_H,B)
$$
for any $\theta \in 
D(\A')$. 
So $\theta \in D(\A')_{<|\A'|-|\A^H|}$ is 
contained in 
$D(\A)$. 
In particular, if there is $\theta \in D(\A')_{|\A'|-|\A^H|} \setminus D(\A)$, then 
for any $\varphi \in D(\A')$, there is $f \in S$ such that 
$\varphi -
f\theta \in D(\A)$. Moreover,
$$
D(\A')_{<|\A'|-|\A^H|} \subset D(\A).
$$
We call this $B$ a \textbf{polynomial $B$ with respect to 
$(\A',H)$}.
\label{B}
\end{theorem}

For the proof of 
Theorem 
\ref{B}, see Proposition 4.41 in \cite{OT}. 
By using several techniques including Theorem \ref{B}, we can show 
the addition-deletion theorems for free arrangements.

\begin{theorem}[Addition-deletion theorem, \cite{T1}]
Let $H \in \A,\ \A':=\A \setminus \{H\}$. Then two of the 
following (1), (2) and (3) imply the third:
\begin{itemize}
    \item [(1)] $\A$ is free with $\exp(\A)=(d_1,\ldots,d_{\ell-1},d_\ell)$.
    \item [(2)] $\A'$ is free with $\exp(\A')=(d_1,\ldots,d_{\ell-1},d_\ell-1)$.
    \item [(3)] $\A^H$ is free with $\exp(\A^H)=(d_1,\ldots,d_{\ell-1})$.
\end{itemize}
Moreover, all the three above hold true if 
$\A$ and $\A'$ are both free.
\label{adddel}
\end{theorem}

It was proved in \cite{A5} that the addition theorem is combinatorial as follows:

\begin{theorem}[Theorem 1.4, \cite{A5}]
Let $\A \ni H$ and $\A':=\A \setminus \{H\}$. If $\A'$ is free, then the freeness of $\A$ depends only on $L(\A)$. 
\label{addcombin}
\end{theorem}

The following is the famous Saito's criterion for the freeness.

\begin{theorem}[Saito's criterion, \cite{Sa}]
Let $\theta_1,\ldots,\theta_\ell \in D(\A)$ be homogeneous derivations. Then 
$\A$ is free with basis $\theta_1,\ldots,\theta_\ell$ if and only if 
they are $S$-independent, and $\sum_{i=1}^\ell \deg \theta_i=|\A|$.
\label{Saito}
\end{theorem}

The following is an easy application of Saito's criterion.

\begin{prop}
Let $\A=\A' \cup \{H\}$ be an arrangement in $\K^3$. 

(1)\,\,
Assume that $\A'$ is free with basis 
$\theta_E,\theta_2,\theta_3$. If $\theta_2 \in D(\A)$, then $D(\A)$ is free with basis 
$\theta_E,\theta_2,\alpha_H\theta_3$.

(2)\,\,
Assume that $\A$ is free with basis 
$\theta_E,\theta_2,\theta_3$. If $\theta_2/\alpha_H \in \Der S$, then $D(\A')$ is free with basis 
$\theta_E,\theta_2/\alpha_H,\theta_3$.

\label{3saito}
\end{prop}

The following results and definition are the key results to show the main results in this article.

\begin{define}[Definition 1.1, \cite{A5}]
We say that $\A$ is \textbf{SPOG (strongly plus-one generated)} with 
exponents $\POexp(\A)=(1,d_2,\ldots,d_\ell)$ and 
\textbf{level} $d$ if there is a minimal free resolution 
$$
0 \rightarrow S[-d-1] \rightarrow \oplus_{i=1}^\ell S[-d_i] \oplus S[-d] \rightarrow D(\A) 
\rightarrow 0.$$
The generators for SPOG arrangements are called the \textbf{SPOG generator}, and among them the 
element whose degree is the level is called the \textbf{level element}.
\label{SPOGdefine}
\end{define}

\begin{theorem}[Theorem 1.4, \cite{A5}]
Let $\A$ be free with exponents $(1,d_2,\ldots,d_\ell)$. If 
$\A':=\A \setminus \{H\}$ is not free, then it is 
SPOG with $\POexp(\A')=(1,d_2,\ldots,d_\ell)$ and 
level $|\A'|-|\A^H|$. 

More 
precisely, if $\theta_E,\theta_2,\ldots,\theta_\ell$ form 
a basis for $D(\A)$, then there is a level 
element $\varphi \in D(\A')_{|\A'|-|\A^H|}$ 
such that $\theta_E,\theta_2,
\ldots,\theta_\ell,\varphi$ form an SPOG genetarot for $D(\A')$.
\label{SPOG}
\end{theorem}

\begin{theorem}[Theorem 5.5, \cite{A5}]
Let $\A'$ be free with exponents $(1,d_2,d_3)$ in $\K^3$. If 
$\A:=\A' \cup \{H\}$ is not free, then it is 
SPOG with $\POexp(\A)=(1,d_2+1,d_3+1)$ and 
level $|\A^H|-1$.

More precisely, if $\theta_E,\theta_2,
\theta_3$ form 
a basis for $D(\A')$, then there 
is a level element $\varphi \in D(\A)_{|\A^H|-1}$
such that $\theta_E, \alpha_H \theta_2,\alpha_H 
\theta_3,\varphi$ form an SPOG generator 
for $D(\A)$.
\label{SPOG2}
\end{theorem}

Now let us recall some definitions and results from \cite{A10}.

\begin{define}[Definition 3.1, \cite{A10}]
For an arrangement $\A$, let $g(\A)$ denote 
the cardinality of 
a minimal set of generators for $D(\A)$. Clearly it is independent of the choice of the set of minimal 
generators. 
\label{ga}
\end{define} 

Moreover we define the following integer for a free arrangement.

\begin{define}[Definition 3.2, \cite{A10}]
Let $\A=\A' \cup \{H\}$ and assume that $\A'$ is free. Let 
$FB(\A')$ be the set of all the 
homogeneous basis for $D(\A')$ and for each 
$B:=\{\theta_1,\ldots,\theta_\ell\} \in FB(\A')$ define 
$$
NT(B,H):=|\{i\mid 1\le i \le \ell,\ 
\theta_i \not \in D(\A)\}|,
$$
and define 
$$
SNT(\A',H):=\min\{NT(B) \mid B \in FB(\A')\}.
$$
\label{MNT}
\end{define}

By using Definition \ref{MNT}, some theory to analyze the logarithmic derivation modules close to the free arrangement is established in \cite{A10} as follows:

\begin{prop}[Proposition 3.4, \cite{A10}]
Let $\A':=\A \setminus \{H\}$ be free with $SNT(\A',H)=s$. Let 
$\theta_E,\theta_2,\ldots,\theta_\ell$ form a basis for $D(\A')$ such that 
$\theta_i \not \in D(\A)\ (2 \le i \le s+1)$ and $\theta_i  
\in D(\A)\ (i \ge s+2)$. Then $g(\A) \ge \ell+s-1$.
\label{genlowbdd}
\end{prop}

\begin{prop}[Proposition 3.5, \cite{A10}]
Let $H \in \A$, $\A':=\A \setminus \{H\}$, $i<j$, and let 
$\A'$ be free with $SNT(\A')=2$. Let $\theta_E,\theta_2,\ldots,\theta_\ell$ be a basis for 
$D(\A')$ such that $\theta_k \in D(\A)\ (k \neq i,j)$. 
Let $\theta_i(\alpha_H) =f_i\alpha_H+g_i B,\ 
\theta_j(\alpha_H)=f_j \theta_H+g_j B$ by 
Theorem \ref{B}.  Then $(g_i,g_j)=1$ and 
$D(\A)$ is generated by 
$\{\theta_k\}_{k \neq i,j} \cup\{\alpha_H \theta_i,\alpha_H \theta_j,g_j\theta_i-g_i\theta_j\}$. In particular, 
$\A$ is SPOG with $\POexp(\A)=
(1,d_2,\ldots,d_{i-1},d_i+1,d_{i+1},\ldots,d_{j-1},d_j+1,d_{j+1},
\ldots,d_\ell)_\le$ for some $i<j$ and level $
d_i+d_j-|\A'|+|\A^H|
$.
\label{fund2}
\end{prop}

\begin{prop}[Proposition 3.8, \cite{A10}]
Let $\A$ be SPOG, $H \in \A$ and $\A':=
\A \setminus \{H\}$. 
If $\A'$ is free, then there is 
SPOG generators $\theta_1=\theta_E,\theta_2,\ldots,\theta_\ell$, a level element  $\varphi$ and two distinct integers 
$1 < s < t \le \ell$ such that 
$$
\theta_E,\theta_2,\ldots,\theta_{s-1},
\theta_s/\alpha_H,\theta_{s+1},\ldots,
\theta_{t-1},
\theta_t/\alpha_H,\theta_{t+1},\ldots,\theta_\ell
$$
form a free basis for $D(\A')$.
\label{fund}
\end{prop}

\section{Main theorem}

Now we have prepared for the proofs of main results.
\medskip

\noindent
\textbf{Proof of Theorem \ref{freepath}}. If $H_1=H_2$ then there is nothing to show, so we may assume that they are distinct. Assume that $\A_1$ is not free and let $\exp(\A)=(1,d_2,\ldots,d_\ell)_\le$.
Then by Theorem \ref{SPOG}, $\A_1$ is SPOG with 
$\POexp(\A_1)=(1,d_2,\ldots,d_\ell)_\le$
and level $|\A|-1-|\A^{H_1}|$.
Then by Proposition \ref{fund}, we know that $\exp(\B)=(1,d_2,\ldots,d_{s-1},d_s-1,d_{s+1},\ldots,d_{t-1},d_t-1,d_{t+1},
\ldots,d_\ell)_\le$ for some $s<t$. Moreover, since $g(\A_1)=\ell+1$,  Proposition \ref{genlowbdd} shows 
that $SNT(\B,H_2)=2$, and the degrees of the basis for $D(\B)$ that is not tangent to $H_2$
are $d_s-1$ and $d_t-1$. So Proposition \ref{fund2} shows that the level of $\A_1$ is $d_s-1+d_t-1-|\A|+2+|\A_1^{H_2}|$. Thus comparing the level,
we have 
\begin{eqnarray*}
d_s-1+d_t-1-|\A|+2+|\A_1^{H_2}|&=&
|\A|-1-|\A^{H_1}|\\
\iff
d_s+d_t-|\A|+|\A_1^{H_2}|&=&
|\A|-1-|\A^{H_1}|\\
\end{eqnarray*}
Let $B$ be the polynomial $B$ with respect to $(\A_1,H_1)$ and let 
$B'$ be that with respect to $(\B,H_2)$. Then by definition and the above equality, we know that 
$$
d_s+d_t-2-\deg B'=\deg B.
$$
Now, let 
$\theta_E,\theta_2,\ldots,\theta_\ell$ be a basis for $D(\B)$ such that 
$\deg \theta_i=d_i\ (i \neq s,t)$ and $\deg \theta_i=d_i-1$ for $i=s,t$. 
By Theorem \ref{B} and the fact that $\theta_s, \theta_t \not \in D(\A_1)$, it holds that $\deg B' \le d_s-1$.
Assume that $\deg B'=d_s-1$. Since $\theta_s 
\not \in D(\A_1)$ and $\deg \theta_s=\deg B'$, $\theta_s(\alpha_{H_2})\equiv B'$ modulo $\alpha_{H_2}$. 
Since $\theta_t$ is not tangent to $H_2$, by putting $\theta_t(\alpha_{H_2})=f_2 B'+\alpha_{H_2}g_2$, 
by Theorem \ref{B}, we may replace $\theta_t$ by $\theta_t-f_2\theta_s$ for some constant $c$ which makes $SNT(\B,H_2)=1$, a contradiction. Thus 
$\deg B' < d_s-1$. So the above equality shows that 
$$
\deg B\ge d_s+d_t-2-(d_s-2)=d_t \ge d_s.
$$
Recall that, by Proposition \ref{fund2}, 
$\theta_E,\theta_2,\ldots,\alpha_{H_2}\theta_s,\ldots,\alpha_{H_2} \theta_t,\ldots,\theta_\ell$ together with a level 
element $\varphi \in S\theta_s+S\theta_t$ form an SPOG generator for $D(\A_1)$. If $\deg B >d_t\ge d_s$, 
then both $\alpha_{H_2}\theta_s,\alpha_{H_2}\theta_t \in D(\A)$ by Theorem \ref{B}. So Theorem \ref{SPOG} shows that at least one derivation in the basis for $D(\A)$ is divisible by $\alpha_{H_2}$, so $\A_2$ is free by the deletion theorem. So we may assume that $\deg B=d_t$. If $d_t>d_s$, then Theorem \ref{B} shows that $\alpha_{H_2}\theta_s \in D(\A)$, implying that $\theta_s \in D(\A_2)$.  Thus $SNT(\B,H_1)=1$, implying that $\A_2$ is free. Hence we may assume that $\deg B=d_s=d_t$. Note that $\deg B=\deg \varphi$ by Theorem \ref{SPOG}, and at least one of $\alpha_{H_2}\theta_s,\alpha_{H_2}\theta_t,\varphi$ 
is not tangent to $H_1$. If it is either $\alpha_{H_2}\theta_s$ or $\alpha_{H_2}\theta_t$, then by Theorem \ref{B}, at least one derivation of the form $\alpha_{H_2} \theta$ form a part of basis by Theorem \ref{SPOG}, hence $\A_2$ is free. Hence let us assume that $\varphi$ is not tangent to $H_2$.
Then $\alpha_{H_2}\theta_s-c_s \varphi$ and $\alpha_{H_2}\theta_t-c_t \varphi$ is tangent to $H_1$ for some constants $c_s,c_t$ by Theorem \ref{B}, and together with $\theta_E$ they form a basis for $D(\A)$ by Theorem \ref{SPOG}. Then we may choose one derivation in the basis as 
$$
c_t(\alpha_{H_2}\theta_s-c_s \varphi)-c_s(\alpha_{H_2}\theta_t-c_t \varphi)
$$
which is divisible by $\alpha_{H_2}$, again confirming that $\A_2$ is 
free. \owari
\medskip


\noindent
\textbf{Proof of 
Theorem \ref{main3}}.
By Theorem \ref{freepath}, the statement is clear if one of $\A_i,\B_j$ is free. So 
we may assume that they are all not free. 
By Theorem \ref{SPOG} and Proposition \ref{fund}, they are all SPOG. 
Let $\theta_E,\theta_1,\theta_2$ be a basis for $D(\A)$ with $1 \le \deg \theta_1 \le \deg \theta_2$, $\theta$ be the level element of $D(\A_1)$. Similarly, 
let $\theta_E,\varphi_1,\varphi_2$ be a basis for $D(\B)$ with $1 \le \deg \varphi_1 \le \deg \varphi_2$, and let $\varphi$ be the level element of $D(\B_3)$. 

Let 
$\exp(\A)=(1,a,b)_\le$. Since $D(\A) \subset D(\B)$, the property of $S$-graded modules with Theorem \ref{Saito} shows that 
$\exp(\B)$ is either 
$(1,a-3,b),\ (1,a-2,b-1),\ (1,a-1,b-2)$ or $(1,a,b-3)$. 

\textbf{Case I}. 
Assume that $\exp(\B)=(1,a-3,b)_<$. Then by Theorem \ref{SPOG2}, 
$\B_3$ is SPOG with generators $\theta_E,\alpha_3\varphi_1,\alpha_3 \varphi_2,\varphi$. Here let $\alpha_i:=\alpha_{H_i}$ for $i=1,2,3$. 
Since $D(\A) \subset D(\B_3)$ and $\deg \alpha_3 \varphi_1=a-2<a<b+1=\deg \alpha_3 \varphi_2$,
$\theta_1=\alpha\alpha_3 \varphi_1$ for some linear form $\alpha$.
Since $\theta_1$ is tangent to $H_1,H_2,H_3$, $\varphi_1$ is tangent to at least one of $H_1$ and $H_2$. Thus Proposition \ref{3saito} (1) shows that either $\B_1$ or $\B_2$ is free, a contradiction.

\textbf{Case II}. 
Assume that $\exp(\B)=(1,a-2,b-1)_<$ and $a<b$. 
Then by Theorem \ref{SPOG2}, 
$\B_3$ is SPOG with generators $\theta_E,\alpha_3\varphi_1,\alpha_3 \varphi_2,\varphi$, so $\POexp(\B_3)=(1,a-1,b)$. First assume that $a<b$. 
Since $D(\A) \subset D(\B_3)$, by the reason of degrees, 
$\theta_1=\alpha \alpha_{H_3}\varphi_1 $ for some linear form $\alpha$. Since $\theta_1 \in D(\A)$, $\varphi_1$ is tangent to at least one of $H_1$ and $H_2$. Then by Theorem \ref{3saito}, $\B_i$ is free for either $i=1,2$, a contradiction. So we may assume $\exp(\A)=(1,a,a)$ and $\exp(\B)=(1,a-2,a-1)$, which is investigated in Case V.

\textbf{Case III}.
Assume that $\exp(\B)=(1,a-1,b-2)$.
First assume that $a<b-1$. Then $\POexp(\B_3)=(1,a,b-1)_<$. Since $D(\A) \subset D(\B_3)$, again $\theta_1=\alpha \alpha_{H_3} \varphi_1$ and the same proof as in Case I and II completes the proof.
So the rest case is $a=b$ or $b=a+1$. The former is investigated in 
Case V.

Hence assume that $b=a+1$. 
Then $\exp(\A)=(1,a,a+1)$ and $\exp(\B)=(1,a-1,a-1)$. By Theorems \ref{SPOG} and \ref{SPOG2}, $\A_1$ is SPOG with $\POexp(\A_1)=(1,a,a+1)$ and level $|\A|-1-|\A^{H_1}|$, and $\B_3$ is SPOG with $\POexp(\B_3)=(1,a,a)$ and level $|\B_3^{H_3}|-1$. 
Since $D(\A_1) \subset D(\B_3)$, $\theta_1 \in D(\A_1)$ is expressed as a linear combination of $\alpha_3 \varphi_1$ and $\alpha_3 \varphi_2$ unless $|\B_3^{H_3}|-1=a$, and in this case $\alpha_3 \mid \theta_1$. So Proposition \ref{3saito} (2) shows that $\A_3$ is free. 
Hence we may assume that $|\B_3^{H_3}|=a+1$, and 
$$
\theta_1=c_1\alpha_3 \varphi_1+
c_2\alpha_3 \varphi_2+\varphi.
$$
Then replace $\varphi$ by $c_1\alpha_3 \varphi_1+
c_2\alpha_3 \varphi_2+\varphi$ to assume that $\theta_1=\varphi$. 
Let $\theta_2=\beta_1\alpha_3\varphi_1+\beta_2\alpha_3\varphi_2+
\beta_3 \varphi$ for linear forms $\beta_i$. Replace $\theta_2$ by $\theta_2-\beta_3 \theta_1$ to express that 
$\theta_2=\beta_1\alpha_3\varphi_1+\beta_2\alpha_3\varphi_2$. 
So Proposition \ref{3saito} (2) shows that 
$\A_3$ is free, a contradiction.

\textbf{Case IV}. 
Assume that $\exp(\B)=(1,a,b-3)$.
First assume that $a<b-3$.
Then we may assume that $\theta_1 =\varphi_1$ since $D(\A) \subset D(\B)$. So $\varphi_1$ is tangent to all $H_i\ (i=1,2,3)$, thus Proposition \ref{3saito} shows that $\B_i$ are free, a contradiction.

Next assume that $b=a+3$. Then 
$\theta_1=c_1 \varphi_1+c_2 \varphi_2$ for some constants $c_i$. Thus by replacing one of the basis for $D(\B)$ by $\theta_1$, Proposition \ref{3saito} shows that $\B_i$ is free for each $i$, a 
contradiction. 

So let us assume that $b=a+2$ and $\exp(\B)=(1,a,a-1)$. Here note that $\deg \varphi_1=a-1$ and $\deg \varphi_2=a$.
If $\theta_1=\alpha \varphi_1$, then $\varphi_1$ is tangent to at least one of $H_1,H_2,H_3$, say $H_1$. So $\B_1$ is free with basis 
$\theta_E,\varphi_1,\alpha_1 \varphi_2$, a 
contradiction. So we may assume that 
$\theta_1=\alpha\varphi_1+c \varphi_2$ for some constant $c\neq 0$. 
Then replacing $\varphi_2$ by $c\varphi_2+\alpha \varphi_1$, 
we may assume that $\varphi_2$ is tangent to all $H_1,H_2,H_3$ and the same argument as the above shows a contradiction.

Finally assume that $b=a+1$ and $\exp(\B)=(1,a,a-2)$.
Here note that $\deg \varphi_1=a-2$ and $\deg \varphi_2=a$.
If $\theta_1=f \varphi_1$, then $\varphi_1$ is tangent to at least one of $H_1,H_2,H_3$ say $H_i$. Then Proposition \ref{3saito} shows that $\B_i$ is free, a contradiction. So we may 
assume that $\theta_1=f\varphi_1+c\varphi_2$ for some non-zero constant $c$. Then the completely same proof 
as when $b=a+2$ completes the proof.

\textbf{Case V}. 
Assume that $\exp(\A)=(1,a,a)$ and $\exp(\B)=(1,a-2,a-1)$.
$D(\A_1)$ is SPOG 
by 
Theorem \ref{SPOG} with exponents $(1,a,a)$ and level $|\A|-1-|\A^{H_1}|$. 
Also $D(\B_3)$ is SPOG by Theorem 
\ref{SPOG2} with exponents $(1,a-1,a)$ and level $|\B_3^{H_3}|-1$. 
If $\theta_1=\alpha \alpha_{H_3}\varphi_1$ or 
$\theta_1=\alpha \alpha_{H_3} \varphi_1+c \alpha_{H_3} \varphi_2$, then 
$\alpha_{H_3} \mid \theta_1$ 
and 
$\A_3$ is free, a contradiction. So it has to hold that $|\B_3^{H_3}|=a+1$, i.e., $\deg \varphi=a$.
By the above, $\theta_i=\beta_i \alpha_3 \varphi_1+c_i \alpha_3 \varphi_2+k_i\varphi_i$ for linear forms $\beta_i$ and some constants $c_i,k_i$ with $k_1 \neq 0,k_2\neq 0$. So $0 \neq k_2\theta_1-k_1\theta_2$ is divisible by $\alpha_3$, thus Proposition \ref{3saito} (2) shows that $\A_3$ is free, a contradiction.\owari
\medskip

If $|\A|-|\B| \ge 4$, then Theorem \ref{main3} does not 
hold true, see the following example.

\begin{example}
Let $\A$ be an arrangement in $\R^3$ consisting of the cone of 
all edges and diagonals of 
the regular 
pentagon. This is known to be free with $\exp(\A)=(1,5,5)$. Since 
$|\A^H|=5 $ for any $H \in \A$, $\A \setminus \{H\}$ 
is not free by Terao's deletion 
theorem.
Let $d\A$ be the deconing of $\A$ with $z=0$, i.e., 
the affine arrangement (a finite set of lines in $\R^2$ that may not 
contain the origin) consisting of 
all edges and diagonals of 
the regular 
pentagon. 
Let $p_1$ be a vertex of 
the regular pentagon. Take all the lines containing 
$p_1$, which we call $\B_1$. 
Pick one edge $H_1 \in \B_1$, and let $p_2 \neq p_1$ be the other 
vertex of the pentagon on $H_1$. Let $p_3$ be the other vertex but $p_1$ next to 
$p_2$ on the pentagon. Then we can choose the edge $H_2$ connecting $p_2$ and $p_3$, and 
the diagonal $H_3$ 
containing $p_1$ and $p_3$. 
Finally, let $H_4$ be the diagonal parallel to $H_1$ and 
through $p_3$. Let $d\B$ consist of $\B_1\ \cup 
\{H_2,H_3,H_4\}$ and $\B$ its coning consisting of 
$7$-planes, which is free with $\exp(\B)=(1,3,3)$. 
Let us show that there are no free paths 
between $\A$ and $\B$. 
If there is, then $\A \setminus \{H\}$ has to be free at least one $H \in \A$, but 
it is impossible since $|\A^H|=5$ for all $H \in \A$, thus $\exp(\A^H)=(1,4) \not 
\subset \exp(\A)=(1,5,5)$ combined with Theorem \ref{adddel}. So Theorem \ref{main3} confirms that 
there are no free arrangements $\C$ with $\B \subsetneq \C \subsetneq \A$.
\label{pentagon}
\end{example}

Finally let us give an application on the free path related to the deformation of the 
Weyl arrangements of rank two. First let us show the following easy corollary.

\begin{cor}
Let $\ell = 3$ and $\A \supset \B$ be both free. 
If the freeness of $\B$ depends on $L(\B)$ and $|\A \setminus \B| \le 3$, then 
so is the freeness of $\A$.
\label{3combin}
\end{cor}

\noindent
\textbf{Proof}.
Clear by Theorems \ref{main3} and \ref{addcombin}. \owari
\medskip

Now let us concentrate our interest on $V=\R^3$. Let $W$ be a Weyl group of rank two, i.e., it is either 
of the type $A_2,\ B_2$ or $G_2$. Let $\Phi^+$ be the corresponding positive system. Then the corresponding 
\textbf{Weyl arrangement} $\A_W$ is defined by 
$$
\A_W=\{\alpha=0\mid \alpha \in \Phi^+ \}.
$$
For example, if $W$ is of the type $A_2$, then 
after a change of coordinates, $\A_W$ is defined by 
$$
x_1x_2(x_1-x_2)=0.
$$
For an integer $k \ge 0$, let us define the 
\textbf{Catalan arrangement} $\Cat^k_W$ by 
$$
\Cat_W^k=\{\alpha=jz\mid \alpha \in \Phi^+,\ -k \le j \le k \} \cup \{z=0\}
$$
and the \textbf{Shi arrangement} $\Shi^k_W$ by 
$$
\Shi_W^k=\{\alpha=jz\mid \alpha \in \Phi^+,\ -k+1 \le j \le k \} \cup \{z=0\}.
$$
Here for the Shi arrangements, $k \ge 1$. Also let $\Phi^+=\Phi_l^+ \cup \Phi_s^+$ be the decomposition of 
the 
roots into the long and short roots. Then define 
\begin{eqnarray*}
\Cat_W^{k_1,k_2}&=&\{\alpha=jz\mid \alpha \in \Phi^+_l,\ -k_1 \le j \le k_1 \} \\
&\cup&
\{\alpha=jz\mid \alpha \in \Phi^+_s,\ -k_2 \le j \le k_2 \}
\cup \{z=0\},
\end{eqnarray*}
\begin{eqnarray*}
\Shi_W^{k_1,k_2}&=&\{\alpha=jz\mid \alpha \in \Phi^+_l,\ -k_1+1 \le j \le k_1 \} \\
&\cup&
\{\alpha=jz\mid \alpha \in \Phi^+_s,\ -k_2+1 \le j \le k_2 \}
\cup \{z=0\},
\end{eqnarray*}
\begin{eqnarray*}
\Cat^{k_1}_W\Shi_W^{k_2}&=&\{\alpha=jz\mid \alpha \in \Phi^+_l,\ -k_1 \le j \le k_1 \} \\
&\cup&
\{\alpha=jz\mid \alpha \in \Phi^+_s,\ -k_2+1 \le j \le k_2 \}
\cup \{z=0\},
\end{eqnarray*}
and
\begin{eqnarray*}
\Shi^{k_1}_W\Cat_W^{k_2}&=&\{\alpha=jz\mid \alpha \in \Phi^+_l,\ -k_1+1 \le j \le k_1 \} \\
&\cup&
\{\alpha=jz\mid \alpha \in \Phi^+_s,\ -k_2 \le j \le k_2 \}
\cup \{z=0\}.
\end{eqnarray*}
They are all proved to be free in \cite{AT2}. Let us prove the following.

\begin{theorem}
(1)\,\,
Let $W$ be of the type $A_2$. Then 
the freeness of $\A$ is combinatorial for all $\Shi^k_W \subset \A \subset \Cat^k_W$ 
or 
$\Cat^k_W \subset \A \subset \Shi^{k+1}_W$
for all $k$.

(2)\,\,
Let $W$ be of the type $B_2$ or $G_2$. Then 
the freeness of $\A$ is combinatorial for all $\Shi^{k_1,k_2}_W \subset \A \subset \Cat^{k_1}_W\Shi^{k_2}_W$,
$\Shi^{k_1,k_2}_W \subset \A \subset \Shi^{k_1}_W\Cat^{k_2}_W$, 
$\Cat^{k_1,k_2}_W \subset \A \subset \Cat^{k_1}_W\Shi^{k_2+1}_W$,  or 
$\Cat^{k_1,k_2}_W \subset \A \subset \Shi^{k_1+1}_W\Cat^{k_2}_W$
for all $k$.
\label{ATmore}
\end{theorem}

\noindent
\textbf{Proof}.
When $W=A_2$, this is already known, e.g. by the results by Yoshinaga's criterion in \cite{Y2}
and Wakamiko's determination of the exponents of the multiarrangements of the type 
$A_2$ in \cite{W}. However, in this case, since $|\Cat_W^k \setminus 
\Shi^k_W|=|\Shi^{k+1}_W \setminus \Cat^k_W|=3$, Theorems \ref{main3}
and \ref{addcombin} also confirms the combinatorial freeness of them.

So let us show the rest cases. 
However in any cases, $\A$ is between two free arrangements $\A_1 \subset 
\A \subset \A_2$ with $|\A_2 \setminus \A_1|\le 3$. Hence Theorem \ref{main3} completes the proof. \owari

\end{document}